\documentclass[twocolumn]{autart}
\usepackage{bm}
\usepackage{array}
\usepackage{amsmath}
\usepackage{graphicx}
\usepackage{amssymb}
\usepackage{arydshln}
\usepackage{graphicx,color,psfrag}
\usepackage{dsfont}
\usepackage[active]{srcltx}
\usepackage{fancyhdr}
\usepackage{enumerate}
\usepackage{multirow}
\usepackage{verbatim}
\usepackage{cite}
\usepackage[usenames,dvipsnames]{pstricks}
\usepackage{epsfig}
\usepackage{upgreek}
\usepackage{multirow}
\usepackage{fancyvrb}
\usepackage{fontenc}
\usepackage{caption}
\usepackage{subcaption}
\usepackage{tikz}
\usetikzlibrary{shapes.geometric, arrows}
\newcommand\script[1]{{\fontfamily{pzc}\fontshape{it}\selectfont#1}}
\VerbatimFootnotes
\definecolor{blau}{rgb}{0,0,0.5}

\newcommand{\R}{\mathbb{R}}
\newcommand{\Prb}{\boldsymbol{\Pr}}
\newcommand{\E}{\boldsymbol{\mathrm{E}}}
\newcommand{\Pm}{\mathcal{P}}
\newcommand{\Tm}{\mathcal{T}}

\newcommand{\tr}{\mathrm{tr}}
\newcommand{\vvec}{\mathrm{vec}}
\DeclareMathOperator*{\argmin}{arg\,min}

\newcolumntype{C}[1]{>{\centering\let\newline\\\arraybackslash\hspace{0pt}}m{#1}}

\begin{document}

\begin{frontmatter}

\title{Co-design of jump estimators and transmission policies for wireless multi-hop networks with fading channels}


\author[prnc]{D. Dolz}\ead{ddolz@uji.es},
\author[quev]{D. E. Quevedo}\ead{dquevedo@ieee.org},
\author[prnc]{I. Pe\~{n}arrocha}\ead{ipenarro@uji.es},
\author[leon]{A. S. Leong}\ead{asleong@unimelb.edu.au},
\author[prnc]{R. Sanchis}\ead{rsanchis@uji.es}
\address[prnc]{Department of Industrial System Engineering and Design, Universitat Jaume I, Castell\'on, Spain}
\address[quev]{Faculty of Electrical Engineering and Information Technology (EIM-E), Paderborn University, Paderborn, Germany}
\address[leon]{Department of Electrical and Electronic Engineering, University of Melbourne, Parkville, Vic. 3010, Australia}



\begin{abstract}
We study transmission power budget minimization of battery-powered nodes in a remote state estimation problem over multi-hop wireless networks. Communication links between nodes are subject to fading, thereby generating random dropouts. Relay nodes help to transmit measurements from distributed sensors to an estimator node. Hopping through each relay node introduces a unit delay. Motivated by the need for estimators with low computational and implementation cost, we propose a jump estimator whose modes depend on a Markovian parameter that describes measurement transmission outcomes over a finite interval. It is well known that transmission power helps to increase the reliability of measurement transmissions, at the expense of reducing the life-time of the nodes' battery. Motivated by this, we derive a tractable iterative procedure, based on semi-definite programming, to design a finite set of filter gains, and associated power control laws to minimize the energy budget while guaranteeing an estimation performance level. This procedure allows us to tradeoff the complexity of the filter implementation with performance and energy use.
\end{abstract}

\end{frontmatter}

\section{Introduction}

Wireless communication technologies have considerably improved in recent years in terms of reliability and transmission rates. This has favored their use for control and estimation purposes~\cite{chen2011guest}, replacing traditional wired technologies. Wireless networks offer several advantages in contrast to wired ones, such as ease of manoeuvre, low cost and self-power. However, wireless links are subject to channel fading that may lead to time-varying delays and packet dropouts~\cite{hespanha2007survey}. These  network-induced issues must be taken into account when designing networked control systems.

Considering remote estimation over networks, Kalman filter approaches may yield optimal performance (for linear dynamical systems), but possibly at the expense of notable implementation and computational complexity, because the time-varying filter gain needs to be computed at each instant in real-time (e.g.~\cite{Sinopoli2004,Schenato08}). Motivated by offering low cost alternatives, we explore in this work the use of precalculated gains that alleviate the computing requirements~\cite{Smith2003,Han2013,Dolz2014b}. As background to our current work, the authors in~\cite{Smith2003} proposed a jump linear estimator whose gains depend on the history of measurement transmission outcomes. Estimation performance improvement is achieved at the expense of increasing estimator complexity (i.e., storage demands and gain selection mechanism). Intermediate complexity approaches were presented in~\cite{Han2013,Dolz2014b}. In~\cite{Han2013} the authors employed a gain dependency on the possible arrival instants and delays for packetized measurements in a finite set, while in~\cite{Dolz2014b} we also considered the multi-sensor case.

Recently, a great deal of attention has focused on designing both the estimator and the transmission policy (this could be regarded as a co-design problem~\cite{smarra2012optimal}) exploiting the relationship between transmission power and dropouts~\cite{Quevedo2012,Gatsis2014}. Higher transmission power leads to lower dropout probabilities, which improves estimation performance. On the other hand, battery life-time is of great importance as battery replacement can be difficult and expensive, and transmission is often the most power consuming task~\cite{estrin2002wireless}. Some relevant works on this topic include~\cite{Shi2012,Ren2014,nourian2014optimal}. These works present methodologies to minimize the estimation error using a Kalman filter while limiting the energy use. These works considered point-to-point communication, i.e., only two nodes (sensor and estimator) are concerned in the data communication. However, e.g., due to the large distance between transmitters and receivers or to the presence of obstacles in the path, simple point-to-point transmission through wireless fading channels may be highly unreliable or extremely power consuming~\cite{Molisch2010}.

With the aim of improving measurement delivery and reducing power budget, in this work we focus on multi-hop wireless networks (see~\cite{Molisch2010}) where some nodes (usually called relays) consciously help to transmit the information from the source to the final destination. These topologies are based on the fact that node data broadcasts are more likely to be acquired from nearby nodes.
In our recent articles ~\cite{leoque13,Quevedo2014} we studied the estimation (using a Kalman filter) and power control problem through multi-hop fading networks. However, these works neglected any transmission delay when hopping through intermediate nodes. Delays can lead to performance issues when estimating the states of a fast dynamical system (for instance, communicating between nodes through IEEE 802.15.4 networks will typically take around 10 ms~\cite{hart}). Motivated by this fact and in the spirit of~\cite{Shi2011}, in the current work we will assume that hopping through each relay introduces an additional unit delay on the data. While in~\cite{Shi2011} the authors presented a two-hop network with two power levels (direct transmission or transmission through relay), here we analyze more general network topologies where multiple relay nodes are present, leading to multiple communication paths between sensors and the estimator node (with different end-to-end delays) and transmission success probabilities.

In this paper we study the transmission power budget minimization of wireless self-powered nodes in a remote state estimation problem for multi-sensor systems over multi-hop networks. Wireless links are subject to
fading leading to random dropouts; hopping through each intermediate node introduces an additional unit delay. We describe this via a finite measurement outcome parameter taken as a finite Markov chain and, based on the network average behavior, we propose a jump linear filter structure. As a difference w.r.t.~\cite{Smith2003}, we use convex optimization in the design of the filter, which allows one to include constraints to fix the number of gains of the jump filter and leads to a trade-off between filter complexity (implementation and computational burden) and estimation performance. We characterize this compromise and give some insights on how to reduce the filter complexity via Lagrange multipliers. We study the co-design problem of minimizing the power budget while guaranteeing a prescribed estimation performance. Since this optimization is non-linear, we derive a greedy algorithm that solves iteratively semi-definite programming problems in order to obtain the set of filter gains and the power transmission laws. 

\textbf{Notation} :  Let $\R$ and $\R_{\geq 0}$ denote respectively the real and  positive real numbers. For matrices $A$ and $B$, $A\preceq B$ means that the matrix $A-B$ is negative semidefinite, and $A\prec B$ means that the matrix $A-B$ is negative definite. The direct sum is represented as $\bigoplus$, thus $A\bigoplus B$ is a block diagonal matrix with $A$ and $B$ on its diagonal. For a given finite set $\mathcal{L}$, $|\mathcal{L}|$ denotes its cardinality. Expected value and probability are denoted as $\E\{\cdot\}$ and $\Prb\{\cdot\}$. The operators $\bigvee$, $\bigwedge$ and $\neg$ represent respectively the logical ``or'', ``and'' and ``not''.

\section{Remote estimation over a multi-hop network}\label{sec:prbapr}
We consider a linear time invariant discrete-time system defined by:
\begin{align}\label{estados}
&x[k+1]=A\,x[k]+B\,w[k],\\
&y_s[k]=c_s\,x[k]+v_s[k],\label{eq:y_cs}
\end{align}
where $x[k]\in\mathbb{R}^n$ is the system state, $y_s[k]\in\mathbb{R}$ is the $s$-th measured output ($s=1,\ldots,n_y$), $w[k]\in\mathbb{R}^{n_{w}}$ is the state disturbance assumed to be a Gaussian signal of zero mean and (known) covariance $\E\{w[k]\,w[k]^T\}=W$, and $v_s[k]\in\mathbb{R}$ is the $s$-th sensor's measurement noise considered as an independent zero mean Gaussian signal with (known) variance $\E\{{v_{s}[k]}^2\}=\sigma_s^2$. For further reference, we define $y[k]\triangleq[y_1[k]\,\cdots\,y_{n_y}[k]]^T$. Also, we assume the pair $(A,C)$ to be detectable, where $C=[c_1^T\,\cdots\,c_{n_y}^T]^T$.

In this work, we study the remote estimation of the system states~\eqref{estados} where the received measurements at the estimator node arrive through an unreliable multi-hop wireless network with fading channels and known topology. We assume that multiple sensors sample the system outputs synchronously and send them independently through the network to a centralized estimator. In the interest of simplicity, we assume that nodes work in a half-duplex mode with mutually orthogonal wireless links, i.e., nodes cannot send and receive at the same time and there is no interference between them~\cite{Molisch2010}. Moreover, we suppose that the nodes access the communication channels with a time division multiple access (TDMA) method using a predefined protocol. Thus, we assume that nodes are time-driven and synchronized. Here, we consider multi-hop wireless networks that can be described via an acyclic directed graph\footnote{Connections between nodes have a direction and each node on a directed path is visited only once.}. This kind of routing has been extensively studied
in the literature, e.g.~\cite{dana2006capacity,alur2011compositional} and is supported by IEEE 802.15.4 transceivers. We denote the set of network nodes by $\displaystyle{\mathcal{N}=\{N_1,\ldots,N_M,N_{M+1}\}}$ with $M>n_y$ being the number of transmitter nodes. $N_1$ to $N_{n_y}$ are the sensor nodes, $N_{n_y+1}$ to $N_M$ are the relay nodes and $N_{M+1}$ refers to the estimator node. While relay nodes are used to retransmit data, sensors can only send their own samples. The network topology is classified and ordered by layers depending on the maximum number of hops (longest path) for a transmission to arrive at the estimator from each node. We assume that the number of different layers is bounded by $\bar{d}+2$ and thus, the maximum number of hops is $\bar{d}+1$. The set of nodes in the $d$-layer is denoted by ${\mathcal{N}_d\triangleq\{N_a\in\mathcal{N}_d:|(N_a,N_{M+1})|=d\}\subset\mathcal{N}}$ where $|(N_a,N_{M+1})|$ stands for the maximum number of hops from $N_a$ to the estimator node. The $0$-layer contains only the estimator node, the $(\bar{d}+1)$-layer includes only sensor nodes, and all other layers may comprise either relay nodes (intermediate nodes that help to retransmit the data) or sensors. 

At each instant $k$, a set of nodes (that transmit in different frequency bands) aggregate all their available measurements in a single time-stamped packet and broadcast it once (without retransmissions) at the same time. Only nodes within a lower layer will accept the transmission (i.e, from $d_1$-layer to $d_2$-layer with $d_1>d_2$), establishing wireless links. The rest of the nodes in the same or higher layers ignore the reception. Thus, a node may receive multiple measurement packets from higher layer nodes and may forward this information to various lower layer nodes. We denote the entire set of wireless links as $\mathcal{I}$, and a single link as $(N_a,N_l)\in\mathcal{I}$. When the dedicated transmission time slot is over, the following set of nodes starts to transmit. After all nodes have attempted to communicate (and before the sampling period has passed), the estimator uses all the received information at instant $k$ to run the state estimation algorithm to be presented in Section~\ref{sec:MJO}. While each sensor transmits the current sampled output, each relay node transmits at instant $k$ only the acquired data at $k-1$.  Whenever a relay node has nothing to retransmit, it frees its channel.


Different from~\cite{Quevedo2013} and as in~\cite{Shi2011}, the transmission protocol implies that communicating through each relay layer, introduces an additional unit delay (equivalent to one sampling period). Direct transmissions to the estimator node do not incurs delays. Thus, a measurement being transmitted at time $k$ by sensor node $N_{s}\in\mathcal{N}_{d+1}$ may arrive at the estimator node with an end-to-end delay of up to $d$ time steps, depending on the number of intermediate layers  visited. The estimator node discards measurements already received.

\begin{figure}[h]
\begin{center}
  \includegraphics[width=\linewidth]{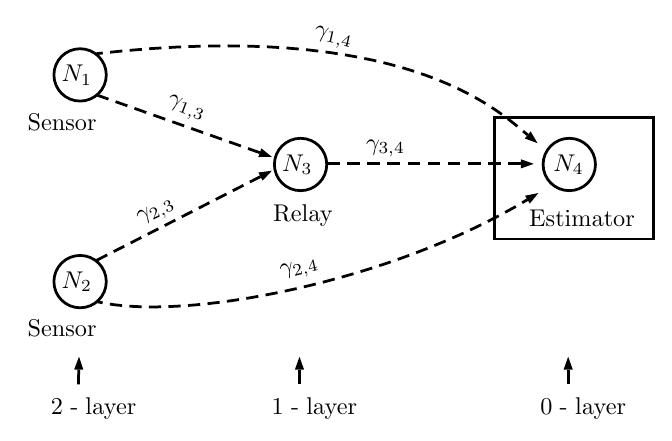}
\caption{Multi-hop network.
}\label{fig_ej_multihop}
\end{center}
\end{figure}

\begin{exmp}\label{ej:multihop_topology} In the multi-hop wireless network of Fig.~\ref{fig_ej_multihop}, 
a measurement leaving $N_1$ may reach the estimator directly (one hop) or may hop twice if it is transmitted through $N_3$. Thus, 
$\bar{d}=1$, with $\mathcal{N}_2=\{N_1, N_2\}$ (direct transmission), $\mathcal{N}_1=\{N_3\}$ (transmission through one relay) and  $\mathcal{N}_{0}=\{N_4\}$ (estimator). If a measurement $y_1[k]$ arrives directly at $N_4$ (at time $k$) and also through $N_3$ (at time $k+1$), then the delayed measurement will be discarded.
\end{exmp}

\section{Transmission outcome model}\label{sec:trans_out}

To model the unreliable transmission through the available wireless links $(N_a,N_l)\in\mathcal{I}$, we introduce the following binary variable:
\begin{equation}\label{def:gamma}
\gamma_{a,l}[k]=\begin{cases}
   1 & \text{if $N_l$ receives a packet from $N_a$ at $k$,} \\
   0 & \text{otherwise.}
  \end{cases}
\end{equation}
Throughout the first part of this work, we assume that each $\gamma_{a,l}[k]$ is an  i.i.d. stochastic process (the effects of fading and power control will be explored in Section~\ref{sec:power}).  The probability of successfully acquiring a transmitted packet is given by:
\begin{equation}\label{def:prb_gamma}
\beta_{a,l}\triangleq\Prb\{\gamma_{a,l}[k]=1\},\quad a,l\in\{1,\ldots,\}.
\end{equation}

We denote by $\tau_s[k]\in\mathbb{N}$ the delay experienced by the $k$-th measurement from sensor $s$ (induced by  hopping in the network) when  accepted at the estimator node. Thus, $\tau_s[k]=d$ means that $y_s[k]$ is accepted, i.e., for the first time received by the estimator  at time $k+d$. The instance $\tau_s[k]>d$ states that the measurement may still be acquired with an induced delay greater than $d$. Since the number of hops is bounded by $\bar{d}+1$, the maximum possible end-to-end delay is $\bar{d}$, i.e., $\tau_s[k]\in\{0,1,\ldots,\bar{d}\}$. Thus, $\tau_s[k]>\bar{d}$ means that $y_s[k]$ is lost.

Let us use $\Gamma_{s,d}^k$ to enumerate the Boolean combinations (logical ``and'' and ``or'' operations) of variables $\gamma_{a,l}$ that define the possible paths a measurement from sensor $N_s$ sent at time $k$ may take to reach and be accepted by the estimator node with a given delay $d$, i.e., with $\tau_s[k]=d$. The possible node-to-node transmission outcomes leading to $\tau_s[k]>d$ are denoted by \reflectbox{\script{$\Gamma$}}$_{s,d}^k$ and can be obtained  by the negation of the disjunction of the corresponding $\Gamma_{s,d}^k$, i.e., \reflectbox{\script{$\Gamma$}}$_{s,d}^k=\neg\left(\bigvee_{\delta=0}^d\Gamma_{s,\delta}^k\right)$.


\begin{exmp} Taking into account the network topology in Fig.~\ref{fig_ej_multihop}, Table~\ref{tab:cases} shows the different outcomes of the node-to-node transmissions leading to each possible end-to-end delay of a measurement sent at time $k$ from sensor $N_1$ (similar applies from sensor $N_2$). 
\begin{table}[ht]
\begin{center}

\caption{Transmission outcomes from $N_1$  in Fig.~\ref{fig_ej_multihop}. '$\times$' denotes indifference occurrence of 1 or 0.}\label{tab:cases}
\vspace{4pt}
\begin{tabular}{|c|c|c|c|cc|}
\hline
From $N_1$& $\Gamma_{1,0}^k$ & $\Gamma_{1,1}^k$ & \reflectbox{\script{$\Gamma$}}$_{1,0}^k$ & \multicolumn{2}{c|}{\reflectbox{\script{$\Gamma$}}$_{1,1}^k$ }\\
\hline
$\gamma_{1,4}[k]$ & 1 & 0 & 0 &0 & 0\\
\hline
$\gamma_{1,3}[k]$ & $\times$ & 1 & $\times$ &0 & $\times$ \\
\hline
$\gamma_{3,4}[k+1]$ & $\times$ & 1 & $\times$ & $\times$ & 0\\
\hline
\end{tabular}
\end{center}
\end{table}
\end{exmp}
Considering the network model described above, the available information at the estimator node at time $k$ are the pairs $(m_{s,d}[k],\alpha_{s,d}[k])$ for all  $s=1,\ldots,n_y$ and $d=0,\ldots,\bar{d}$, where
\begin{equation}\label{ec:m_s_ds}
m_{s,d}[k]=\alpha_{s,d}[k]\,y_s[k-d]
\end{equation}
and $\alpha_{s,d}[k]$ is a binary variable such that
\begin{equation}\label{def:alpha_d_i}
\alpha_{s,d}[k]=\begin{cases}
   1 & \text{if } y_{s}[k-d] \text{ is received at time } k, \\
   0 & \text{otherwise.}
  \end{cases}
 \end{equation}
When $\alpha_{s,d}[k]=1$, the measurement sent at time $k-d$ from sensor $N_s$ has experienced a delay of $\tau_s[k-d]=d$. If $y_s[k-d]$ has not yet arrived at time $k$, then $m_{s,d}[k]=0$. Since delayed copies (already received measurements with a higher delay) are discarded, $\alpha_{s,d}[k]$ is equal to zero if $\alpha_{s,d-\delta}[k-\delta]=1$ for some integer $\delta\in\{1,\ldots,d\}$.

Let us now introduce a vector $\theta_{s,d}[k]$ which models the successful reception of $y_s[k-d]$  during the interval $\{k-{d},k-{d}+1,\ldots,k\}$:
\begin{equation}
\theta_{s,d}[k]=\begin{bmatrix}\alpha_{s,0}[k-d] & \alpha_{s,1}[k-d+1] & \cdots & \alpha_{s,d}[k]\end{bmatrix}.\label{eq:Markov_theta3}
\end{equation}
In the following section, we will study how to design a Markovian jump filter for the proposed network scenario. Our estimator will allow us to explore trade-offs between estimation performance and estimator complexity.

\section{Markovian jump filter}\label{sec:MJO}

To take into account the reception of delayed measurements up to $\bar{d}$ , we propose to use an aggregated model
\begin{equation}\label{estadosaug}
\bar{x}[k+1]=\bar{A}\bar{x}[k]+\bar{B}w[k],
\end{equation}
with $\bar{x}[k]=\begin{bmatrix}x[k]^T&\cdots&x[k-\bar d]^T\end{bmatrix}^T$ and $(\bar{A},\bar{B})$ appropriate augmented matrices.

\subsection{Building the Markov chain}
We denote by $\bar{\tau}\in\mathbb{N}$ a parameter that allows to adjust the length to look back in time. If $\bar{\tau}=\bar{d}$ we consider all the possible delayed measurements. If $\bar{\tau}<\bar{d}$ we narrow the historical interval to just take into account acquired measurements with a delay lower than $\bar{d}$, while if $\bar{\tau}>\bar{d}$ we allow to look further back in time, even if no measurements with a higher delay than $\bar{d}$ will be received. Then, concatenating vectors $\theta_{s,d}[k]$ in~\eqref{eq:Markov_theta3} for $d=0,\ldots,\bar{\tau}$ and for all sensors, we have
\begin{equation}\label{eq:Markov_theta}
\theta[k]=\begin{bmatrix}\theta_{1,0}[k]&\cdots&\theta_{1,\bar{d}}[k]&\theta_{2,0}[k]&\cdots&\theta_{n_y,\bar{\tau}}[k]\end{bmatrix}
\end{equation}
where $\theta[k]$ is a binary column vector, of length $n_\theta=\frac{(\bar{\tau}+1)(\bar{\tau}+2)}{2}n_y$, representing the full set of measurements successfully received from $k-\bar{\tau}$ to $k$ and fulfilling the conditions
\begin{equation} \label{eq:cond_Theta}
\sum_{\delta=\bar{d}+1}^{\bar{\tau}}\alpha_{s,\delta}[t-\delta+d]=0,\quad \|\theta_{s,d}[k]\|_1\leq1,
\end{equation}
where
${\|\theta_{s,d}[k]\|_1=\sum_{\delta=0}^{\bar{\tau}}\alpha_{s,\delta}[t-\delta+d]}$.
These conditions describe the fact that measurements only arrive with a delay up to $\bar{d}$ and that delayed copies are discarded. The possible occurrences of $\theta[k]$  lie within a finite set, i.e., ${\theta}[k]\in\Theta=\{\vartheta_0,\vartheta_1,\ldots,\vartheta_{r}\}$ with\footnote{$\lceil\cdot\rceil$
is the operator that rounds its argument to the nearest
positive integer (including zero) towards infinity.}
\begin{equation}\label{eq:Deltak}
r={\left((\bar{d}+2)!(\bar{d}+2)^{\left\lceil\bar{\tau}-\bar{d}\right\rceil}\right)}^{n_y}-1.
\end{equation}
Each $\vartheta_i$ (for $i=0,\ldots,r$) represents one of the possible combinations of the historical measurement transmission outcomes. 
From the definition of $\theta[k]$, and taking into account the fact that the occurrences of $\alpha_{s,d}[k]$ depend on node-to-node transmission outcomes (defined by $\gamma_{a,l}[k]$, and which are assumed i.i.d., see~\eqref{def:prb_gamma}), we can conclude that $\theta[k]$ is a homogeneous Markov process (see~\cite{bremau99}). We shall assume that it is ergodic.

\begin{rem}
$\theta[k]$ captures the measurement transmission outcomes at times $\{k-\bar{\tau},\ldots,k\}$ where measurements with a delay up to $\bar{d}$ are expected to be acquired. If the probability of any measurement from sensor $N_s$ to be accepted by the estimator after $d+1$ hops (with $d\leq\bar{d}$) is zero, i.e., $\Prb\{\Gamma_{s,d}^k\}=0$, the Markov chain $\theta[k]$ is not irreducible because the states containing the information about measurements with delay $d$ can never be reached from any other state. In this case, we can perform a state space reduction (by removing the infesible states) leading to ${\theta[k]^\prime\in\Theta^\prime\subset\Theta}$ with $|\Theta^\prime|<r$, where $\Theta^\prime$ contains only the reachable states. Thus, $\theta^\prime[k]$ is irreducible. To ensure ergodicity we must guarantee that $\theta[k]^\prime$ is aperiodic too. Since $\theta[k]^\prime$ is irreducible we only need one aperiodic state to imply that it is aperiodic~\cite{bremau99}. Let us assume that the estimator is constantly receiving measurements with the same delay that makes $\theta[k]^\prime$ converge to a stationary state, i.e., the period to return to itself is one, which makes this state aperiodic. Therefore, $\theta[k]^\prime$ is ergodic.
 \end{rem}
The following result shows how to obtain the transition probabilities of $\theta[k]$.
\begin{lem}\label{lem:p_ij}
The elements of the transition probability matrix $\Lambda=[p_{i,j}]$ of $\theta[k]$ with $p_{i,j}\triangleq\Prb\{\theta[k]=\vartheta_j\big{|}\theta[k-1]=\vartheta_i\}$ can be computed as:
\begin{align}\label{ec:p_ij_theta}
 \hspace{-7pt} p_{i,j}&=\Prb\left\{\varphi(k,j)\wedge\varphi(k-1,i)\right\}/\Prb\left\{\varphi(k-1,i)\right\}
\end{align}
with
\begin{align}
&\varphi(k,j)=\left(\bigwedge_{\{s,d,\delta\}\in\mathcal{D}_1(\vartheta_j)} \hspace{-5pt}\Gamma_{s,\delta}^{k-d}\right) \bigwedge \left(\bigwedge_{\{s,d\}\in\mathcal{D}_2(\vartheta_j)} \hspace{-5pt}\reflectbox{\script{$\Gamma$}}_{s,d}^{k-d}
\right),\label{ec:varphi}
\end{align}
$s=1,\ldots,n_y$, $d=0,\ldots,\bar{\tau}$, $\delta=0,\ldots,{d}$ and
\begin{align*}
\mathcal{D}_1(\vartheta_j)&\triangleq\{s,d,\delta:\theta[k]=\vartheta_j,\,\delta\leq\bar{d},\,\alpha_{s,\delta}[k-d+\delta]=1\},\\ \mathcal{D}_2(\vartheta_j)&\triangleq\{s,d:\theta[k]=\vartheta_j,\,\|\theta_{s,d}[k]\|_{1}=0\}.
\end{align*}
\end{lem}
\begin{pf}
$\varphi(k,j)$ is a function returning the possible followed paths by the measurements  sent at times $\{k-\bar{\tau},\ldots,k\}$ leading to the state $\theta[k]=\vartheta_j$. It is composed of the conjunction for each path $\Gamma_{s,\delta}^{k-d}$ when there is a measurement with $\tau_{s}[k-d]=\delta$, which is imposed  by $\{s,d,\delta\}\in\mathcal{D}_1(\vartheta_j)$, and those in \reflectbox{\script{$\Gamma$}}$_{s,\delta}^{k-d}$ that imply the fact that $\tau_{s}[k-d]>\delta$, which is imposed  by $\{s,d\}\in\mathcal{D}_2(\vartheta_j)$. Given the above, $\Prb\{\theta[k]=\vartheta_j\}$ can be rewritten in terms of $\Prb\{\varphi(k,j)\}$ and applying the definition of  conditional probability, we obtain~\eqref{ec:p_ij_theta}.
\end{pf}

\begin{rem}
$p_{i,j}$ is calculated using the node-to-node successful transmission probabilities in~\eqref{def:prb_gamma} since $\Gamma_{s,d}^k$ and \reflectbox{\script{$\Gamma$}}$_{s,d}^k$ are determined by Boolean combinations of variables $\gamma_{a,l}$. Note that when $d>\bar{d}$ we have that \reflectbox{\script{$\Gamma$}}$_{s,d}^k=$\reflectbox{\script{$\Gamma$}}$_{s,\bar{d}}^k$.
\end{rem}

Let us now define the measurement availability matrix at time $k$ as:
\begin{equation}\label{def:alpha}
\alpha[k]=\psi(\theta[k])\triangleq\bigoplus_{s=1}^{n_y}\left(\bigoplus_{d=0}^{\bar{d}}\alpha_{s,d}[k]\right).
\end{equation}
The possible values of $\alpha[k]$ belong  to a set of the form $
{\alpha}[k]\in\Xi=\{\eta_0,\eta_1,\ldots,\eta_{q}\}$ with $q=2^{\bar n_y}-1$ and
where $\eta_i$ (for $i=1,\ldots,q$) denotes each possible combination. 

\begin{exmp}\label{ex:p_ij} Consider the network topology in Fig.~\ref{fig_ej_multihop}. Let us show how to compute the transition probability, when $\bar{\tau}=\bar{d}=1$, between the states
$\theta[k-1]=\begin{bmatrix}1&0&0&0&0&1\end{bmatrix}^T \textrm{ and } \theta[k]=\begin{bmatrix}0&1&0&1&0&1\end{bmatrix}^T$. $\theta[k]=\vartheta_4$ models the situation that $y_1[k]$ is not acquired  at time $k$ (i.e., $\tau_1[k]>0$), $y_1[k-1]$ was received at time $k-1$ (i.e., $\tau_1[k-1]=0$), $y_2[k]$ is acquired time at $k$ (i.e., $\tau_2[k]=0$) and $y_2[k-1]$ is received  at time $k$ (i.e., $\tau_2[k-1]=1$). Similarly, $\theta[k-1]=\vartheta_3$ corresponds to the situation that $\tau_1[k-1]=0,\, \tau_1[k-2]>1,\, \tau_2[k-1]>0$, and $\tau_2[k-2]=1$.
Using expression~\eqref{ec:p_ij_theta}, the transition probability is given by
$ p_{3,4}=\frac{\Prb\left\{\varphi(k,4)\wedge\varphi(k-1,3)\right\}} {\Prb\left\{\varphi(k-1,3)\right\}}$ where
$\varphi(k,4)=\reflectbox{\script{$\Gamma$}}_{1,0}^k\wedge\Gamma_{1,0}^{k-1}\wedge\Gamma_{2,0}^k\wedge\Gamma_{2,1}^{k-1}$
and $\varphi(k-1,3)=\Gamma_{1,0}^{k-1}\wedge\reflectbox{\script{$\Gamma$}}_{1,1}^{k-2}\wedge\reflectbox{\script{$\Gamma$}}_{2,0}^{k-1}\wedge\Gamma_{2,1}^{k-2}$.
Finally, considering the paths described in Table~\ref{tab:cases} and the node-to-node successful transmission probabilities~\eqref{def:prb_gamma}, we obtain ${p_{3,4}=\Prb\{\neg\gamma_{1,4}[k]\wedge\gamma_{2,4}[k]\wedge\gamma_{2,3}[k-1]\wedge\gamma_{3,4}[k]\}}$ ${=(1-\beta_{1,4})\,\beta_{2,4}\,\beta_{2,3}\,\beta_{3,4}}$, see~\eqref{def:gamma}.
\end{exmp}

Finally, using vector $\bar{m}_s[k]$ we express the availability of the measurements at time $k$ from sensor $s$ sent from $k-\bar{d}$ to $k$ as $ \bar{m}_{{s}}[k]=[m_{s,0}[k]\,\cdots\,m_{s,\bar{d}}[k]]^T$,
where $m_{s,d}[k]$ is as defined in~\eqref{ec:m_s_ds}. Using $\alpha[k]$, we rewrite the received measurement information at time $k$ as:
\begin{equation}\label{ec:bar_m_k}
\bar{m}[k]={\alpha}[k]\left({\bar{C}}\bar{x}[k]+\bar{v}[k]\right)
\end{equation}
where $\bar{m}[k]=[{\bar m_1[k]}^T\;\cdots\;{\bar m_{n_y}[k]}^T]^T$, $\bar{v}[k]=[\bar v_{1}[k]\;\cdots\;\bar v_{\bar n_y}[k]]^T$ with $\bar v_{s}[k]=[v_{s}[k]\;\cdots\;v_{s}[k-\bar{d}]]$, and ${\bar{c}}_{s}=[{\bar{c}}_{s,0}\,\cdots\,{\bar{c}}_{s,{\bar{d}}}]^T$ with  ${\bar{c}}_{s,d}=[0_{1\times n\cdot d}\; c_{s}\;0_{1\times n\cdot(\bar d-d)}]^T$ (and $c_s$ as defined in~\eqref{eq:y_cs}) are the rows of matrix ${\bar{C}}$. In~\eqref{ec:bar_m_k}, $\bar{v}[k]$ is the measurement noise vector with covariance $\E\{\bar{v}[k]\bar{v}[k]^T\}=V=\bigoplus_{s=1}^{n_y}\left(\bigoplus_{{d}=0}^{\bar d}\sigma_{s}^2\right).$

\subsection{Proposed filter}
We propose the state estimation algorithm
 \begin{align}
\hat{\bar{x}}[k]&=\bar A\hat{\bar{x}}[k-1]+L[k](\bar m[k]-\alpha[k]\bar C\bar A\hat{\bar{x}}[k-1]).\label{ec:estimator_update}
\end{align}
When no measurement is available, the estimator is run in open loop. Otherwise, the state estimation is updated with the updating gain matrix $L[k]$. Considering~\eqref{estadosaug} and~\eqref{ec:bar_m_k}-\eqref{ec:estimator_update}, the dynamics of the estimation error, defined as
$\tilde{x}[k]=\bar{x}[k]-\hat{\bar{x}}[k]$, is
\begin{align}\label{ec:estimator_error}
\tilde{x}[k]=&(I-L[k]\alpha[k]\bar C)\left(\bar{A}\tilde{x}[k-1]+\bar{B}w[k-1]\right)\nonumber\\
&-L[k]\alpha[k] \bar{v}[k].
\end{align}
One of the aims of this work is to compute gain matrices $L[k]$ to obtain acceptable estimation performance while requiring low computing and storage capabilities.
Using predefined gain filters~\cite{Smith2003,Dolz2014b} instead of Kalman filters~\cite{Schenato08} helps to alleviate the on-line computational burden, and also allows for dealing with e.g. uncertain systems. The authors in~\cite{Smith2003} illustrated that the Kalman filter gains depend on the history of sensor availability. In the current work we extend results in~\cite{Smith2003} to include multisensor transmission and delayed measurements. Thus, we propose to adapt the gains to ${\theta}[k]$ such as:
\begin{equation}\label{def:Lt}
L({\theta}[k])=\begin{cases}
   0  &\text{if } \psi({\theta}[k])=0,\\
      L_i   &\text{if } {\theta}[k]=\vartheta_i,\,\psi(\vartheta_i)\neq 0.
  \end{cases}
\end{equation}
We will compute the gain matrices off-line leading to the finite set
\begin{equation}\label{setL}
L({\theta}[k])\in\mathcal{L}=\{L_{0},\ldots,L_{r}\}.
\end{equation}

\begin{rem}\label{rem:rprime}  We measure the filter complexity by the number of different stored gains in $\mathcal{L}$, i.e., by $|\mathcal{L}|$ . 
Notice from (\ref{def:Lt}) that $L(\theta[k])=0$ is repeated $r_0=(r+1)/(\bar{d}+2)^{n_y}-1$ times. Thus,  $|\mathcal{L}|$ is at most $|\mathcal{L}|=r+1-r_0$. For fixed values of $\bar{\tau}$ and $\bar{d}$, we can further reduce the filter complexity by imposing equality constraints over the set $\mathcal{L}$, e.g. by setting ${L}_{i}={L}_{j}$ for some $i\neq j$. This will be explored in the sequel.
\end{rem}

\section{Filter design}\label{sec:obsvdsgn}

To design the filter, we first note that the Markov chain $\{{\theta}[k]\}$  has a stationary distribution (due to ergodicity) that satisfies $\pi=\pi\Lambda$, where $\pi=[\pi_1,\ldots,\pi_{r}]$ and $\pi_i=\Prb\{{\theta}[k]=\vartheta_i\}$ are the probabilities of being at a given state. Based on this, the next result characterizes the evolution of the state estimation error covariance matrix.

\begin{thm}[\cite{Dolz2014a}]\label{teor:Pt}
Let $P_{j}[k]=\E\{\tilde x[k]\tilde x[k]^T|{\theta}[k]=\vartheta_j\}$ (with $j=1,\ldots,r$) be the modal state estimation error covariance matrix updated at time $k$ with information ${\theta}[k]=\vartheta_j$. We then have:
\begin{align}\label{ec:Pt_j}
P_{j}[k]=&\sum_{i=0}^{r}p_{i,j}\frac{\pi_i}{\pi_j}\left(F_{j}(\bar{A}P_{i}[k-1]\bar{A}^T+\bar{B}W\bar{B}^T)F_{j}^T\right)\nonumber\\
&+\sum_{i=0}^{r}p_{i,j}\frac{\pi_i}{\pi_j}X_{j}VX_{j}^T
\end{align}
with $ F_{j}=I-L_{j}\,\psi(\vartheta_j)\bar C$ and $X_{j}=L_{j}\psi(\vartheta_j)$. Furthermore, the expected value of the state estimation error covariance matrix is given by:
\begin{align}\label{ec:Pt}
\E\{\tilde x[k]\tilde x[k]^T\}=\sum_{j=0}^{r}P_{j}[k]\pi_j.
\end{align}
\end{thm}
The above theorem defines a recursion on the modal covariance matrices in~\eqref{ec:Pt_j}, that we write as ${P_j[k]=\mathfrak{E}_j\{\Pm[k-1]\}}$ where
$\Pm[k]\triangleq(P_{0}[k],\ldots,P_{r}[k])$ and $\mathfrak{E}_{j}\{\cdot\}$ is the linear operator over all the modal covariance matrices that gives~\eqref{ec:Pt_j}. Thus, we write the full recursion (for $j=0,\ldots,r$) as $\Pm[k]=\mathfrak{E}\{\Pm[k-1]\}$ where
$\mathfrak{E}\{\cdot\}\triangleq(\mathfrak{E}_0\{\cdot\},\ldots,\mathfrak{E}_{r}\{\cdot\})$.
To compute a steady state solution, one must address the problem of finding a set of filter gains that satisfy the Riccati equation $\mathfrak{E}\{\Pm[k]\}=\Pm[k]$. When the filter gains are allowed to jump with all the states of the Markov chain, \cite{Smith2003} and~\cite{Han2013} show how to obtain their explicit values. However, the methods in~\cite{Smith2003,Han2013} do not consider cases where the filters share the same gain for different modes of the Markov chain (i.e, ${L}_{i}={L}_{j}$ for some $i\neq j$). Thus,~\cite{Smith2003,Han2013} do not allow us to explore the trade-offs between storage complexity and estimation performance. To overcome this issue, we adopt the following alternative optimization problem:
\begin{equation}\label{problema}
\underset{\mathcal{L},\Pm}{\text{minimize}}\quad \Tm\{\Pm\} \quad \mathrm{s.t.} \quad \mathfrak{E}\{\Pm\}-\Pm\preceq 0
\end{equation}
with $\Pm\triangleq(P_{0},\ldots,P_{r})$, $\Tm\{\Pm\}=\tr\left(\sum_{j=0}^{r}P_{j}\pi_j\right)$ and where $\mathfrak{E}\{\Pm\}-\Pm\preceq 0$ denotes that $\mathfrak{E}_j\{\Pm\}-P_j\preceq 0$ for all $j=1,\ldots,r$. We showed in~\cite{Dolz2014a} that the feasibility of problem~\eqref{problema} is a sufficient condition to guarantee the boundedness of $\E\{\tilde x[k]\tilde x[k]^T\}$. The iteration $\Pm[k]=\mathfrak{E}\{\Pm[k-1]\}$ converges to the unique positive semi-definite solution ${\Pm}$ for the given $\mathcal{L}$, obtained both from~\eqref{problema} (similar arguments can be found in~\cite{costa2006discrete}). Let us now state some necessary conditions which must be satisfied in order for ~\eqref{problema} to have a solution.

\begin{thm}\label{teor:exist}
For problem~\eqref{problema} to have a solution, the transition probabilities of  ${\theta}[k]$ (and thus the node-to-node successful transmission probabilities in~\eqref{def:prb_gamma}) must fulfill the following constraints
\begin{subequations}\label{con_nec}
\begin{align}
p_{jj}\cdot \rho(\bar A)^2-1\leq0,\,\,\,\,\forall j:\psi(\vartheta_j)=0,\label{cond_nec_p0}\\
p_{jj}\cdot \rho(\bar A_i)^2-1<0,\,\,\,\,\forall j:\psi(\vartheta_j)=\eta_i,\,\,i\in N\hspace{-0.07cm}D,\label{cond_nec_pi_nueva}
\end{align}
\end{subequations}
where $p_{ij}$ are the probabilities defined in~\eqref{ec:p_ij_theta}, $\rho(\bar A)$ denotes the spectral radius of matrix $\bar A$, $N\hspace{-0.07cm}D$ is the set of reception scenarios $\eta_i$ from which $(\bar{A},\eta_i\bar{C})$ is non-detectable, and $\rho(\bar{A}_i)$ is the spectral radius of the unobservable subspace of $\bar{A}$ from the reception scenario $\eta_i$\footnote{Note that $\bar A_i=\mathcal{O}^T\bar{A}\mathcal{O}$, where $\mathcal{O}=\mathrm{ker}\left(\begin{bmatrix}(\eta_i\bar{C})^T&(\eta_i\bar{C}\bar{A})^T& \cdots &  (\eta_i\bar{C}\bar{A}^{n-1})^T\end{bmatrix}^T\right)$ is the unobservable subspace.}.
\end{thm}
\begin{pf} A necessary condition for the constraint in~\eqref{problema} to hold is
$p_{j,j}(I-L_{j}\,\psi(\vartheta_j)\bar C)\bar{A}P_{j}\bar{A}^T(I-L_{j}\,\psi(\vartheta_j)\bar C)^T-P_j\preceq0$ for all $j=0,\ldots,r$, see~\eqref{ec:Pt_j}. i) When $\psi(\vartheta_j)=0$, we have that $p_{j,j}\bar{A}P_{j}\bar{A}^T-P_j\preceq0$. Then, a necessary condition for the existence of a solution is~\eqref{cond_nec_p0}. ii) Following the elimination lemma~\cite{Scherer2001}, the existence of matrices $L_j$ and $P_j$ such that the strict inequality holds, is equivalent to the existence of a matrix $P_j$ with $\psi(\vartheta_j)=\eta_i$ such that ${(\eta_i\bar C \bar A)^\bot}^T\left(p_{jj}(\bar A)^TP_j\bar A-P_j\right)(\eta_i\bar C \bar A)^\bot \succ 0$
where $(\eta_i\bar C \bar A)^\bot$ is a basis for the null space of $(\eta_i\bar C \bar A)$, thereby containing a basis of the unobservable subspace from reception scenario $\eta_i$. Thus, a necessary condition for the existence of the filter is~\eqref{cond_nec_pi_nueva}.
\end{pf}

\subsection{Complexity vs performance}\label{subsec:trade}
We can reduce the storage complexity of the filter by reducing the measurement outcome history $\bar{\tau}$ taken into account or by imposing the same filter gain for different states $\theta_k$ at the cost of a worst performance. The aim is to obtain the lower number of different filter gains that lead to a given prescribed filter performance $\Tm\{\Pm\}$. In order to share a filter gain for a couple of states one must set constraints like ${L}_{i}-{L}_{j}=0$ for some $i\neq j$ over $\mathcal{L}$. Thanks to the convergence property $\Pm=\mathfrak{E}\{\Pm\}$, we can rewrite problem~\eqref{problema} with the corresponding equality constraints as
\begin{equation}\label{LiLj}
\underset{\mathcal{L},\Pm}{\text{minimize}}\quad \Tm\left\{\mathfrak{E}\{\Pm\}\right\} \quad \mathrm{s.t.} \quad {L}_{i}={L}_{j},\,i\neq j.
\end{equation}
One can solve the previous optimization problem for all the possible pairings and then choose the shares that guarantee the prescribed performance with the lower number of different gains. This procedure is highly consuming and other approaches must be explored. Here we address the problem  by first solving the unconstrained problem and then look for recursively the new pair that can share the gain without affecting too much the achieved performance. Each possible equality constraint to be added in the problem will not affect at the same level the achieved performance $\Tm\{\Pm\}$, so only one of the pairings will be kept for further recursions while the desired performance is still satisfied.
\begin{rem}
Problems~\eqref{problema} and~\eqref{LiLj} can be rewritten as a linear matrix inequality (LMI) optimization problem that is solved iteratively, e.g., using the  SeDuMi solver~\cite{sturm1999using}. Details are omitted for brevity but can be found in~\cite{Dolz2014a}. In the general case, each iteration has $V_D=n(\bar{d}+1)(n(\bar{d}+1)+1)(r+1)+nn_y(\bar{d}+1)^2|\mathcal{L}|$ decision variables and $r+1$ sparse LMIs. The latter can be combined into a single LMI of size up to $V_M=\left(n(4r+5)+n_y\right)(\bar{d}+1)(r+1)$. Then, the computational complexity of  SeDuMi is in $\mathcal{O}(V_D^{2}V_M^{2.5}+V_M^{3.5})$. Note that not only the storage requirements are reduced when imposing equality constraints over $\mathcal{L}$, but also the off-line computational complexity to obtain a solution is alleviated.
\end{rem}

The previous procedures are highly time consuming, as one must solve a huge number of optimization problems to decide which pairing is the less harmful to filter performance. Here we try to reduce that time by means of analyzing which pairing (constraint) to be added is less harmful for the achieved performance, and we measure this by means of the sensitivity of the w.r.t. the constraint to be added.   For a fixed $\bar{\tau}$ (a simple rule is to set it to $\bar{d}$), let us assume that we have fixed ${L}_{i}-{L}_{j}=0$  (other gain equality constraints could already have been included), and let us quantify the the effect of perturbing this equality over $\Tm\{\Pm\}$.
The method of Lagrange multipliers tackles this problem by defining the Lagrange functional
\begin{equation}
\Lambda(\mathcal{L},\Pm,\lambda)=\Tm\left\{\mathfrak{E}\{\Pm\}\right\} +\vvec(\lambda)^T\vvec({L}_{i}-{L}_{j})
\end{equation}
and solving the homogeneous equations resulting from the partial derivatives w.r.t  $(\mathcal{L},P,\lambda)$, where $\vvec(\cdot)$ is the vectorization operator. It is known that the resulting value of $\lambda$ quantifies the effect of modifying the new constraint $L_i=L_j$ on the achievable cost index~\cite{bazaraa2013nonlinear}. In this case, we obtain:
\begin{equation}\label{eq:lambda}
\lambda=\sum_{q=0}^r\pi_qZ_q\bar{C}^T\left(M_i-M_j\right)+K\sum_{q=0}^r\pi_q\left(Y_{i,q}-Y_{j,q}\right)
\end{equation}
with
\begin{align*}
&K=\sum_{q=0}^rZ_q[k-1]\bar{C}^T\left(M_i+M_j\right)\left(\sum_{q=0}^rY_{j,q}+Y_{i,q}\right)^{-1},\\
&M_i=\sum_{i\in\mathcal{K}}p_{q,i}\psi(\vartheta_{i})^T,\,\,\,Z_q=\bar{A}\bar{P}_q\bar{A}^T+\bar{B}W\bar{B}^T,\\
&Y_{i,q}=\sum_{i\in\mathcal{K}}p_{q,i}\psi(\vartheta_{i})X_q\psi(\vartheta_{i})^T, \,\,\,X_q=\bar{C}Z_q\bar{C}^T+V
\end{align*}
where $\mathcal{K}$ is a set containing the index of the states $\vartheta_i$ that already share $L_i$ (similarly with $L_j$) and $\bar{P}_q$ is the modal covariance of the state $\vartheta_q$ obtained with the additional constraint $L_i=L_j$. As we want to avoid the complete computation of the optimization problem before choosing a pairing, $\bar{P}_q$ is not available. In order to make use of this information just for pairing purposes, we can approximately quantify the Lagrange multipliers using the modal covariance $\bar{P}_q$ obtained in a previous step, free of this constraint. With this, a lower time consuming iterative procedure to decrease the filter complexity is as follows:
 \textit{(i)} Solve initially the problem (unconstrained or with any arbitrarily fixed pairings) and store the covariance $\bar{P}_q$; \textit{(ii)} For all the couple of possible pairings to be added, compute the lagrange multipliers~\eqref{eq:lambda} using the previous matrices $\bar{P}_q$;  \textit{(iii)} Fix the additional gain equality constraint that leads to the lowest value $\|\lambda\|_2$ for further iterations and return to step \textit{(i)}. The algorithm stops when performance $\Tm\{\Pm\}$ in step \textit{(i)} exceeds the prescribed value, and then one must remove the last pairing to guarantee the prescribed performance. The algorithm also stops when the number of different gains is acceptable for implementation issues. This idea is further explored in Section~\ref{sec:ej}.

\section{Transmission over fading channels}\label{sec:power}

So far, we have assumed that the node-to-node successful transmission probabilities in~\eqref{def:prb_gamma}, used to compute the predefined estimator gains, were known and time-invariant. However, in wireless networks with fading channels, the probability of successfully acquiring at $N_l$ a transmitted packet from $N_a$ depends on the fading channel gain $h_{a,l}[k]\in\Omega_{a,l}\subset\R_{\geq 0}$ and on the transmission power $u_a[k]\in[0,\,\bar{u}]$ as per (cf., ~\eqref{def:prb_gamma})
 \begin{equation}\label{def:prb_gamma_hu}
\Prb\{\gamma_{a,l}[k]=1|h_{a,l}[k]=h,\,u_a[k]=u\}=f_{a,l}(h\,u)
\end{equation}
where the function $f_{a,l}$ is monotonically increasing and differentiable, and depends on the modulation scheme employed, e.g., see~\cite{Proakis1995}.

Let us aggregate in vector $H_a[k]\in\Omega_a$ the outgoing fading channel gains from node $N_a$, i.e., $H_a[k]=[h_{a,l}[k]\,\,\cdots\,\,h_{a,m}[k]]^T$ with $\Omega_a=\Omega_{a,l}\times\cdots\times\Omega_{a,m}$  where $\{(N_a,N_l),\ldots,(N_a,N_m)\}\subset\mathcal{I}$. Assuming that each node $N_a$ knows\footnote{This can be attained in practice by means of channel estimation algorithms, see references in~\cite{Quevedo2012,Quevedo2014}.} $H_{a}[k]$. We shall focus on local power control policies of the form
\begin{equation}\label{eq:ugen}
       u_a[k]=\kappa_a(H_{a}[k]),
\end{equation}
where $\kappa_a:\Omega_a\rightarrow[0,\,\bar{u}]$ is a parameterized and integrable function over $\Omega_a$. To use the model in~\eqref{def:prb_gamma} one could seek to control the power to reach the same constant successful transmission probability at each instant (which might not be achieved). Alternatively, $\beta_{a,l}$ can also be treated as the average behavior of the communication channel over an infinity-time window. For further references we denote by $\mathcal{U}_a$ the set that contains the parameters of $\kappa_a(\cdot)$. We shall next show how to retrieve these values.

\subsection{Network average behavior}

The fading in channel $(N_a,N_l)$ is a stochastic process that might be correlated with other channels (representing some spatial correlation~\cite{agrawal2009correlated}) and that we assume to have a time-invariant disitribution. With this, we denote by $g_a(H_a)$ the joint probability density function of $H_a[k]=H_a$ which is considered to be known. Then, we can compute the transmission probability $\beta_{a,l}$ defined in~\eqref{def:prb_gamma} using:
\begin{align}\label{def:prb_gamma_avg}
 \beta_{a,l}&= \int_{\Omega_a}\,g_a(H_a)\,f_{a,l}\left(h_{a,l}\,\kappa_a(H_a)\right)\mathrm{d}H_a.
\end{align}

\begin{rem} Typically, IEEE 802.15.4 transceivers have a finite number of predefined discrete power levels~\cite{Shi2012,Quevedo2010}. One way to implement this using~\eqref{eq:ugen} is to use a piecewise constant policy wherein each of the power levels is used in different ranges of fading gain values.
\end{rem}

\section{Co-design}\label{sec:trade}
\begin{figure}[h]
\begin{center}
 \includegraphics[width=\linewidth]{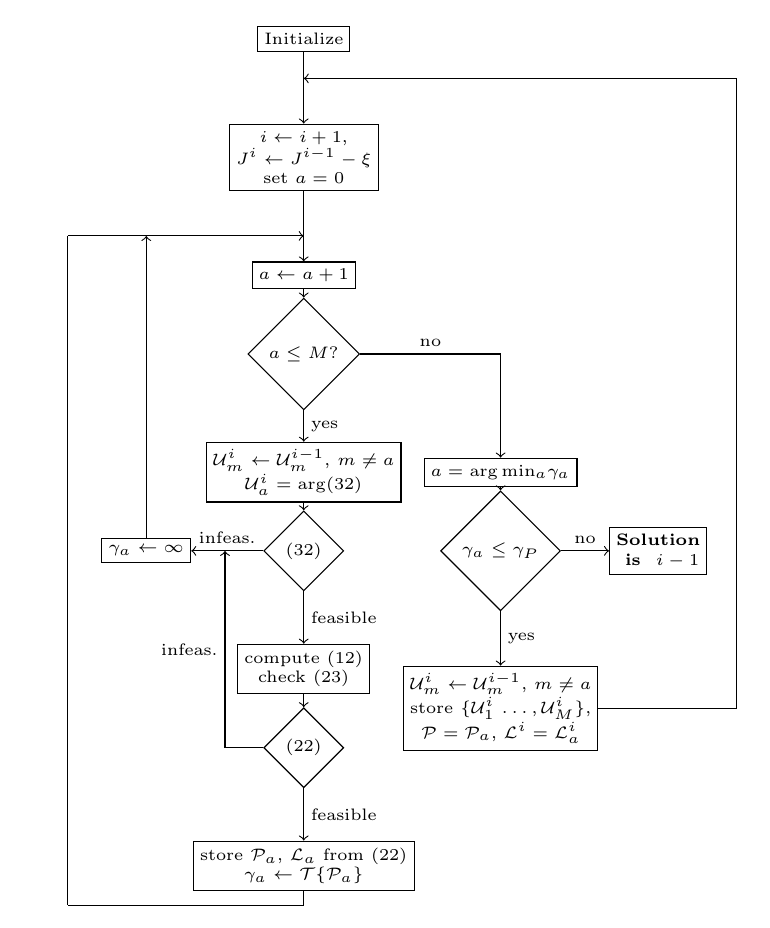}
\caption{Greedy algorithm flow diagram.}\label{fig_algo}
\end{center}
\end{figure}

Transmitter nodes are often self-powered, and thus, preserving battery life is an important concern. Motivated by this, we will next show how to compute off-line the parameters that define the power control policies in~\eqref{eq:ugen} and a set of filter gains guaranteeing the performance of the estimator. Our aim is to guarantee a certain estimation performance $\gamma_P$ so that the network transmission power usage is minimized. Let us characterize the power budget by:
\begin{align}\label{eq:J}
J(U)=\sum_{a=1}^M\mu_a\E\{u_a\}
\end{align}
where $\mu_a\in\R$ and the expected value (average) of the transmission power of each node $N_a$ is:
\begin{equation}\label{eq:U_expt}
\E\{u_a\}=\int_{\Omega_a}g_a(H_a)\kappa_a(H_a)\mathrm{d}H_a.
\end{equation}
Then, the synthesis  problem can be formulated as:
\begin{equation}\label{problemacode}
\begin{aligned}
& \underset{\mathcal{L},\Pm,\mathcal{U}}{\text{minimize}}
& & J(U) \\
& \text{subject to}
& & \Tm\{\Pm\}\leq\gamma_P,\,\,\,\mathfrak{E}\{\Pm\}-\Pm\preceq 0,\,\,\,\eqref{def:prb_gamma_avg},\\
& & & u_a[k]\in[0,\,\bar{u}],\quad \forall a=1,\ldots,M.
\end{aligned}
\end{equation}
This is a nonlinear optimization problem,  as the average node-to-node transmission probabilities~\eqref{def:prb_gamma_avg}, which are employed in deriving the transition probabilities of $\theta[k]$, depend on the power control strategies. This kind of problem can be tackled, for instance, by brute force using a gridding approach, by means of heuristic optimization with genetic algorithms, or by implementing a greedy algorithm. As a fast way to obtain a possibly suboptimal solution, in this work we propose the use of a greedy algorithm. A greedy algorithm is a tree search where at each step only the branch that locally optimizes the problem fulfilling some heuristic is explored, in the hope that this choice will lead to a globally optimal solution. Thus, this kind of algorithm never comes back to previous incumbent solutions to change the search path, and globally optimal solutions are not guaranteed. Nonetheless, often good solutions will be found. The proposed greedy algorithm is as follows:
\begin{description}
\item[Step 1.] Set $i\leftarrow0$. For a given $\bar{u}$, choose $M$ sets of power control parameters $\mathcal{U}_a^0, a=1,\dots,M$ to maximize each power transmission $\E\{u_{a}\}$. For given constants $\mu_a$, define index $J^0\triangleq\sum_{a=1}^M\mu_a\E\{u^0_{a}\}.$ Choose a small positive parameter value $\xi$.

\item[Step 2.] Set $i\leftarrow i+1$ and $J^i\leftarrow J^{i-1}-\xi$. For $a=1$ to $M$ repeat Step 3, then go to Step 4.

\item[Step 3.] Set $\mathcal{U}_m^i\leftarrow\mathcal{U}_m^{i-1}$ for all $m\neq a$ with $m=1,\ldots,M$. Obtain $\mathcal{U}_a^i$ as the power control parameter set resulting from the optimization problem
\begin{equation}\label{problemastep3}
\begin{aligned}
& \underset{\mathcal{U}_a}{\text{maximize}}
& & \prod_{ l\,:\,(N_a,N_l)\in\mathcal{I}}  \beta_{a,l}
\\
& \text{subject to}
& & J^i-\sum_{\substack{m=1\\m\neq a}}^M\mu_m\E\{u_{m}^i\}-\mu_a\E\{u_{a}^i\}=0,\\
& & & u_a[k]\in[0,\ \bar{u}],\quad \forall a=1,\ldots,M.
\end{aligned}
\end{equation}
with $\beta_{a,l}$ as defined in~\eqref{def:prb_gamma_avg}. If this problem has no solution, then set $\gamma_a\leftarrow\infty$. Otherwise, compute the transition probabilities given in~\eqref{ec:p_ij_theta} and check conditions~\eqref{con_nec}. If they are not fulfilled, then set $\gamma_a\leftarrow\infty$. Otherwise solve optimization problem~\eqref{problema} including the corresponding gain equality constraints
leading to the desired $|\mathcal{L}|$. If the problem is infeasible, then set $\gamma_a\leftarrow\infty$. Otherwise, store $\Pm_a$,  $\mathcal{L}_a$ and set $\gamma_a\leftarrow\mathcal{T}\{\Pm_a\}.$

\item[Step 4.] Set $\displaystyle{a\leftarrow\underset{a}\argmin\;\gamma_a}.$ If $\gamma_{a}\leq\gamma_P$, then set $\mathcal{U}_m^i\leftarrow\mathcal{U}_m^{i-1}$ for all $m\neq a$; store $\mathcal{U}^i=\{\mathcal{U}_1^i,\ldots,\mathcal{U}_M^i\}$, $\Pm^i=\Pm_a^i$ and $\mathcal{L}^i=\mathcal{L}_a^i$, and go to Step 2. Otherwise, exit with the best solution found in iteration $i-1$.
\end{description}
The algorithm starts in Step 1 by considering the most favorable power control policy, i.e., where the probabilities of receiving packets are highest (higher transmission powers). Then, at each iteration (Step 2, 3 and 4) it tries to reduce the power budget~\eqref{eq:J} while the feasibility of the problem~\eqref{problemacode} is preserved, see Fig.~\ref{fig_algo}. Each time it can be reduced, in Step 3,  we first translate  the effect of the reduction of the power budget to a single node. This leads to as many power transmission policies (control parameters $\mathcal{U}_a^i$ ) as there are nodes. To obtain the new power control law characterized by $\mathcal{U}_a^i$, we solve~\eqref{problemastep3} where we only modify the transmission policy of $N_a$ such that the successful transmission rate is maximized, whilst fulfilling the new power budget\footnote{The optimization problem~\eqref{problemastep3} can be solved using nonlinear optimization algorithms (e.g., \verb+fmincon+ of MATLAB).}. Then, with the new set of power control parameters $\mathcal{U}_a^i$ for $N_a$ and the already existing $\mathcal{U}_m^i$, the algorithm computes the transition probabilities of the reception scenario model $\theta[k]$ and verifies the fulfillment of the filter existence necessary conditions developed in Theorem~\ref{teor:exist}. If the latter hold, then we design the state estimator by solving~\eqref{problema} including the predefined gain equality constraints. Once this has been done for the $M$ nodes, in Step 4 the proposed heuristic selects the solution with the lowest estimation performance index $\gamma_a$, i.e., the solution that generates a larger future search space. The algorithm ends when all the obtained $\gamma_a$ are higher than the prescribed upperbound $\gamma_P$.

\begin{rem}
The associated Lagrange multiplier of $\mathrm{tr}\left(\sum_{j=0}^{r}P_{j}\pi_j\right)\leq\gamma_P$ is non-negative by definition, which means that lower power budgets can only be obtained by degrading estimation performance.
\end{rem}

\begin{rem} The filter storage complexity is supposed to be a priori fixed in the previous algorithm. To enable the decision on which states of $\theta[k]$ share the same gain, we can modify Step 3 by including an iterative procedure that adds at each step a new equality gain constraint $L_i=L_j$ such that the correspondent Lagrange multiplier~\eqref{eq:lambda} is at a minimum (in the sense of $\|\lambda\|_2$) until the constraint on $|\mathcal{L}|$ is reached. 
\end{rem}

\section{Numerical studies}\label{sec:ej}
\begin{figure}
    \centering
    \begin{subfigure}[b]{0.47\textwidth}
        \centering \includegraphics[width=\linewidth]{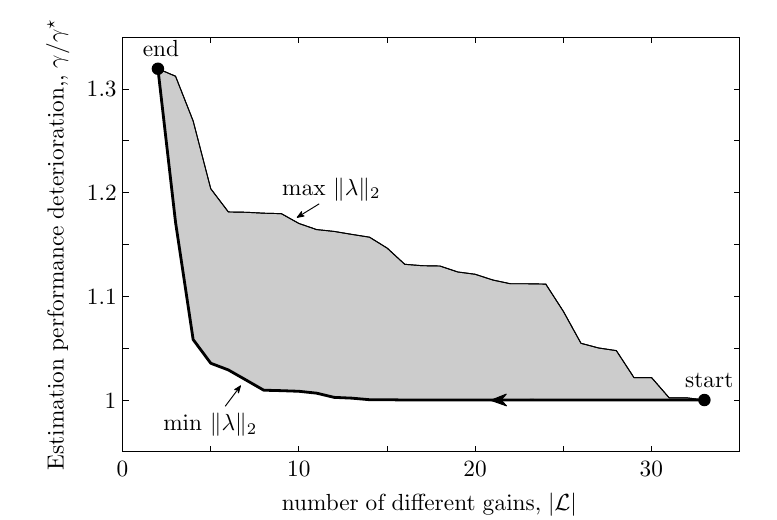}
        \caption{Estimation performance vs. filter complexity, where ${\gamma^\star=0.037}$.}\label{fig:trade}
    \end{subfigure}
    \begin{subfigure}[b]{0.47\textwidth}
        \centering \includegraphics[width=\linewidth]{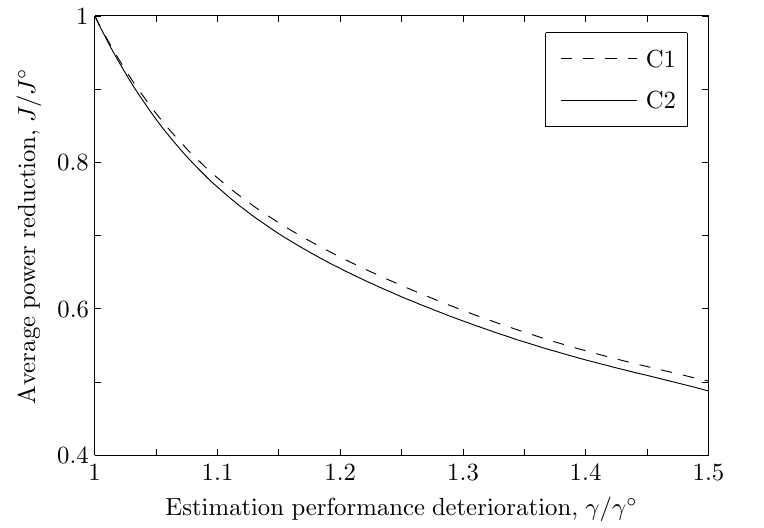}
\caption{Power budget vs. estimation
performance for constant (C1) and saturated inverted channel gain (C2) transmission power policies, where $J^\circ=30$ and $\gamma^\circ=0.023$.}\label{fig:code}
    \end{subfigure}
    \caption{Achievable trade-offs.}\label{fig:ejem8}
\end{figure}

We consider the following (unstable) system
\begin{align*}
A&=\begin{bmatrix}1.05&-0.1\\0.74&1.05\end{bmatrix},\, B=\begin{bmatrix}0.01&0.13\\0.01&0.08\end{bmatrix},\, C=\begin{bmatrix}0.53&0.39\\0.72&0.35\end{bmatrix}.
\end{align*}
The state disturbance and sensor noise covariances are 
\begin{equation*}
W=\begin{bmatrix}
    0.26  & -0.003\\
   -0.003  &  0.25\\
\end{bmatrix},\;
\begin{bmatrix}
\sigma_1^2\\\sigma_2^2
\end{bmatrix}=\begin{bmatrix}
0.0086\\0.0079
\end{bmatrix}.
\end{equation*}
Measurements are acquired through the multi-hop network in Fig.~\ref{fig_ej_multihop} that may induce up to a unit delay in the end-to-end transmission. Thus, the number of states of the Markov chain $\theta[k]$ is $|\Theta|={\left((1+2)!\right)}^2=36$ (see~\eqref{eq:Deltak}). Nodes transmit using BPSK modulation (see~\cite{Quevedo2012}) with $b=4$ bits and a transmission power bounded by $\bar{u}=10$. We consider correlated fading channels with ${h}_{1,4}[t]=(h_{1,3}[t]+h_{3,4}[t])/100$ and ${h}_{2,4}[t]=(h_{2,3}[t]+h_{3,4}[t])/10$ where ${h_{1,3}[t],\,h_{2,3}[t],\,h_{3,4}[t]}$ follow an independent exponential distribution (Rayleigh fading) with means $\bar{h}_{1,3}=1,\, \bar{h}_{2,3}=0.3$ and $\bar{h}_{3,4}=0.5$.
We denote the estimation performance index $\tr\left(\sum_{i=0}^{35}C_x P_{i}C_x^T\pi_i\right)$ as $\gamma$, where $C_x=[I_{n}\,0_{n\times(n\cdot \bar{d})}]$ selects the covariance corresponding to $x[k]-\hat{x}[k|k]$. First, let us assume that each node uses a constant power $u_1[k]=u_2[k]=u_3[k]=5$. Under this scenario, the presence of a relay helps to improve the performance index obtained while solving~\eqref{problema} (without gain equality constraints) from $\gamma=0.112$ to $\gamma^\star=0.037$, where $\gamma^\star$ was retrieved with 33 different gains. In this case, the estimation performance index obtained with the Kalman filter is $\gamma_{\mathrm{Kal}}=0.034$ (where $\gamma_{\mathrm{Kal}}= \tr(C_x\E\{\tilde{x}[k]\tilde{x}[k]^T\}C_x^T)$, which is $7\%$ lower. However, the Kalman filter needs to execute up to 976 floating-point operations at every time instant, while the off-line method only requires at most 64 (independent of the gain grouping strategy). We can further reduce the  filter complexity (storage requirements) using the results in Section~\ref{subsec:trade} where we sequentially impose the gain equality constraint leading to a minimum deterioration of $\gamma$, i.e., the one with the smallest associated Lagrange multiplier in the sense of $\|\lambda\|_2$. Fig.~\ref{fig:trade} explores the achieved trade-offs showing that the proposed decision rule leads to the best estimation performances. Also, we notice that we can reduce the filter complexity from 33 to 14 different gains without affecting the estimation performance index.

Using the greedy algorithm developed in Section~\ref{sec:trade}, we now explore the design trade-offs between power budget and estimation performance when all the filter gains are different. For that purpose, we focus on constant (C1) and saturated inverted channel gain (C2) transmission power policies (see~\cite{Quevedo2012}). Strategy C2 is implemented using the information of the realizations $h_{1,3}[t],\,h_{2,3}[t]$ and $h_{3,4}[t]$ whereas C1 only uses statistical information. Fig.~\ref{fig:code} shows that at the expense of  only a 10\% deterioration of the best estimation performance, we can reduce the power budget by 25\%. Note that, for the given example, power control strategy C2 only reduces the power budget by up to 4\%  with respect to C1 at the expense of requiring the value of the channel fading gain at each instant.

\section{Conclusions}\label{sec:conclu}
In this work we developed  a model for multi-hop networked estimation with fading channels. Random dropouts are generated due to fading links, while delays are introduced while hopping through relay nodes. We introduced a  finite Markovian process that captures the network behavior by keeping track of a finite number of measurement transmission outcomes. Using the average network behavior, we conceived a jump filter whose complexity can be selected as a trade-off between storage requirements and estimation performance. To keep the network operation power-efficient, we designed power policies and filter gains to minimize the power budget while guaranteeing a certain estimation error performance index. The design is carried out with an iterative procedure based on semi-definite programming.  Numerical results suggest that: 1) intermediate relays help to reduce the power budget for prescribed estimation performance, 2) increasing estimator complexity allows in general for power savings, 
and 3) not too many different filter gains are needed to achieve similar estimation performance than Kalman filter approach while offering a much lower computational burden. Future research may include design of topology reconfiguration methods and distributed power control strategies.



\begin{ack}                               
This work has been funded by projects DPI2011-27845-C02-02 from MICINN,  PI15734,  E-2015-15 and P1$\cdot$1B2015-42 from Universitat Jaume I .   
\end{ack}

\bibliographystyle{plain}
\bibliography{bibliografiacov,bib_power}
\end{document}